\input amstex
\magnification=\magstep0
\documentstyle{amsppt}
\pagewidth{6.0in}
\input amstex
\topmatter
\title
Wilf's ``Snake Oil" Method Proves an Identity \\
in The Motzkin Triangle
\endtitle
\author
Tewodros Amdeberhan, Moa Apagodu and Doron Zeilberger
\endauthor
\address Tulane University, Department of Math., 6823 St. Charles Av,
New Orleans, LA 70118 \endaddress
\email tamdeber\@tulane.edu \endemail
\address Virginia Commonwealth University,
Department of Math. \& Applied Math.,
Richmond, VA 23284  \endaddress
\email mapagodu\@vcu.edu \endemail
\address Rutgers University, Department of Math., 110 Frelinghuysen Rd, Piscataway, NJ 08854
\endaddress
\email zeilberg\@math.rutgers.edu \endemail
\abstract We give yet-another illustration of using Herb Wilf's 
{\it Snake Oil Method}, by proving a certain identity between the entries of the so-called 
{\it Motzkin Triangle}, that arose in a recent study of enumeration of certain classes
of integer partitions. We also briefly illustrate how this method can be applied
to general `triangles'.
\endabstract
\dedicatory
Dedicated to the memory of Herb Wilf
\enddedicatory
\endtopmatter
\def\({\left(}
\def\){\right)}

\document
\noindent
Our starting point was a certain conjecture, concerning the so-called 
\it simultaneous core \rm partitions, 
found in a recent preprint [2, Conjecture 11.5]. It reads:
\bigskip
\noindent
\bf Conjecture.  \it Let $s$ and $d$ be two coprime positive integers. Then the number of $(s,s+d,s+2d)$-core partitions is given by
$$\sum_{k=0}^{\lfloor\frac{s}2\rfloor}\binom{s+d-1}{2k+d-1}\binom{2k+d}k\frac1{2k+d}.$$ \rm

\noindent
We then paid particular interest to the special cases $d=1$ (see [1]), resulting in the \it Motzkin numbers \rm proven in [2] and by Yang-Zhong-Zhou [5], and $d=2$ initiating yet another link [1, Problem 11.6] to the Motzkin triangle which we now state. 

\bigskip
\noindent
\bf Problem. \rm The \it Motzkin triangle \rm $T(n,k)$ of numbers is defined according to the rules:

(1) $T(n,0)=1$;

(2) $T(n,k)=0$ if $k<0$ or $k>n$;

(3) $T(n,k)=T(n-1,k-2)+T(n-1,k-1)+T(n-1,k)$.

\smallskip
\noindent
\it Prove the identity (this is sequence $A026940$ in OEIS  [3]))
$$\sum_{k=0}^nT(n,k)T(n,k+1)=\sum_{k=0}^n\binom{2n}{2k+1}\binom{2k+1}k\frac1{k+2}.$$  \rm

Let us first observe that any such identity is nowadays {\it automatically provable},
thanks to the so-called {\it Wilf-Zeilberger algorithmic proof theory}, but
it is still fun to prove it, whenever possible, the old-fashioned way, by purely
{\it human} means. We will do this, by using what Herb Wilf called the Snake-Oil method
\rm [4, Section 4.3].

\pagebreak

\smallskip
\noindent
Recall that the {\it Constant Term} of a Laurent polynomial, $P(x)$, is the
coefficient of $x^0$. For example, $CT[4/x+3+5x]=3$.
\smallskip
\noindent
For motivation, let's look at a few known examples.
\smallskip
\noindent
(1) $\sum_{k=0}^n\binom{n}k\,x^k=(1+x)^n$ and hence 
$\text{CT}\left[\frac{(1+x)^n}{x^k}\right]=\binom{n}k$, 
the famous binomial coefficients as entries in the familiar Pascal's triangle (see $A007318$ in OEIS [3]).

\smallskip
\noindent
(2) $\sum_{k=0}^{n+2}C(n,k)\,x^k=(1+x)^n(1-x)$, 
this is one variant among the Catalan triangles (see the sequences $A008315$ and $A037012$ in OEIS [3]).

\smallskip
\noindent
(3) $\sum_{k=0}^{2n}t(n,k)\,x^k=(1+x+x^2)^n$, the trinomial triangle (see $A027907$ in OEIS [3]).

\bigskip
\noindent
Going back to the Motzkin triangle
, we return to  our Problem 
by first \it extending \rm the definition of the Motzkin triangle from $k=0,1,\dots,n$ to 
$k=0,1,\dots,2n+2$ as a skew-symmetric sequence: 
$$T(n,k)= - T(n,2n-k+2).$$
\bf Note. \rm $T(n,n+1)=0$ and the extended Motzkin triangle 
(we continue to denote by $T(n,k)$) obeys the 
\it same \rm recurrence. As a result, it is easy to construct the generating function
$$\sum_{k=0}^{2n+2}T(n,k)x^k=(1+x+x^2)^n(1-x^2).$$
Or, equivalently, for $k\in\{0,1,\dots,2n+2\}$,
$$T(n,k)=\text{CT}\left(\frac{(1+x+x^2)^n(1-x^2)}{x^k}\right).$$
So, the stage is now set and The Snake Oil method can be brought to bear:
$$\align \sum_{k=0}^nT(n,k)T(n,k+1)&=\frac12\sum_{k=0}^{2n+2}T(n,k)T(n,k+1) \\
&=\frac12\sum_{k=0}^{2n+1}T(n,k)\,\cdot\text{CT}\left(\frac{(1+x+x^2)^n(1-x^2)}{x^{k+1}}\right) \\
&=\frac12\,\text{CT}\left[\frac{(1+x+x^2)^n(1-x^2)}{x}\sum_{k=0}^{2n+1}T(n,k)x^{-k}\right] \\
&=\frac12\,\text{CT}\left[\frac{(1+x+x^2)^n(1-x^2)}{x}\left(1+\frac1x+\frac1{x^2}\right)^n\left(1-\frac1{x^2}\right)\right]  \\
&=\frac12\,\text{CT}\left[\frac{(1+x+x^2)^{2n}(1-x^2)}{x^{2n+1}}\right] -
\frac12\,\text{CT}\left[\frac{(1+x+x^2)^{2n}(1-x^2)}{x^{2n+3}}\right]  \\
&=\frac12\, T(2n,2n+1)-\frac12\, T(2n,2n+3) \\
&=\frac12\, T(2n,2n-1); \endalign$$
where the last equality is due to $T(2n,2n+1)=0$ and $T(2n,2n+3)=-T(2n,2n-1)$. 

\smallskip
\noindent
We pause for a moment to appreciate a striking similarity between the two identities, 
$$\sum_{k=0}^n\binom{n}k\binom{n}{k+1}=\binom{2n}{n+1} \qquad \text{and} \qquad
\sum_{k=0}^nT(n,k)\,T(n,k+1)=\frac12\,T(2n,2n-1),$$
involving coefficients in the Pascal's triangle and the current Motzkin's triangle, respectively.

\pagebreak

\bigskip
\noindent
On the other hand, if we expand $(1+x+x^2)^{2n}=((1+x)+x^2)^{2n}=\sum_{k=0}^{2n}\binom{2n}kx^{2k}(1+x)^{2n-k}$ then by reverse-engineering the expression for $\frac12\,T(2n,2n-1)$ from above, we are led to
$$\align \frac12\,T(2n,2n-1)
&=\frac12\,\text{CT}\left[\frac{(1+x+x^2)^{2n}(1-x^2)}{x^{2n-1}}\right] \\
&=\frac12\,\text{CT}\left[\frac{(1+x+x^2)^{2n}}{x^{2n-1}}\right]
-\frac12\,\text{CT}\left[\frac{(1+x+x^2)^{2n}}{x^{2n-3}}\right] \\
&=\frac12\sum_{k=0}^{n-1}\binom{2n}k\binom{2n-k}{2n-2k-1} 
-\frac12\sum_{k=0}^{n-2}\binom{2n}k\binom{2n-k}{2n-2k-3} \\
&=\frac12\sum_{k=0}^{n-1}\binom{2n}{2k+1}\binom{2k+1}k
-\frac12\sum_{k=0}^{n-2}\binom{2n}{2k+3}\binom{2k+3}k \\
&=n+\frac12\sum_{k=1}^{n-1}\binom{2n}{2k+1}\binom{2k+1}k
-\frac12\sum_{k=1}^{n-1}\binom{2n}{2k+1}\binom{2k+1}{k-1} \\
&=n+\sum_{k=1}^{n-1}\binom{2n}{2k+1}\binom{2k+1}k\frac1{k+2}
=\sum_{k=0}^{n-1}\binom{2n}{2k+1}\binom{2k+1}k\frac1{k+2}, \endalign$$
which is exactly the right-hand side of our problem. This completes the proof. In fact, we have improved the assertion of Problem 11.6 because of our success in evaluating the two sums into the \it simpler \rm form $\frac12T(2n,2n-1)$. Therefore, we may formulate our conclusion as the next result.

\bigskip
\noindent
\bf Theorem 1. \it The following identities hold true:
$$\sum_{k=0}^nT(n,k)T(n,k+1)=\sum_{k=0}^n\binom{2n}{2k+1}\binom{2k+1}k\frac1{k+2}=\frac12\,T(2n,2n-1).$$ \rm

\bigskip
\noindent
A litmus test (or a cannon measure, if you prefer) to the quality of a good technique is perhaps its enlightenment, simplicity and implications. Indeed, in our case, the linear operator CT offers both a clue to and a proof for an \it effortless \rm generalization of Theorem 1. The Motzkin triangle persists!
\bigskip
\noindent
\bf Theorem 2. \it The following identity holds true:
$$\sum_{k=0}^{\lfloor\frac{s}2\rfloor}\binom{s+d-1}{2k+d-1}\binom{2k+d-1}k\frac1{k+d}=\frac1d\,T(s+d-1,s).$$ \rm
\bf Proof. \rm This is completely analogous to what has been demonstrated earlier. To wit,
$$\align \frac1d\,T(s+d-1,s)
&=\frac1d\,\text{CT}\left[\frac{(1+x+x^2)^{s+d-1}}{x^s}\right] 
-\frac1d\,\text{CT}\left[\frac{(1+x+x^2)^{s+d-1}}{x^{s-2}}\right]  \\
&=\frac1d\sum_{k\geq0} \binom{s+d-1}k\binom{s+d-k-1}{s-2k} 
-\frac1d\sum_{k\geq0} \binom{s+d-1}k\binom{s+d-k-1}{s-2k-2} \\
&=\frac1d\sum_{k\geq0} \binom{s+d-1}{2k+d-1}\binom{2k+d-1}k 
-\frac1d\sum_{k\geq0} \binom{s+d-1}{2k+d+1}\binom{2k+d+1}{k} \\
&=\sum_{k\geq0}\binom{s+d-1}{2k+d-1}\binom{2k+d-1}k\frac1{k+d}.
\endalign$$
The proof is complete. $\square$

\pagebreak

\bigskip
\noindent
Is there more? Yes, here is a \it bonus\rm! 
As a nice implication of the preceding results, 
the above Conjecture may be stated much more succinctly.
\bigskip
\noindent
\bf Conjecture. \it If $s, d\geq 1$ are coprime integers, then the number of $(s,s+d,s+2d)$-core partitions equals
$$\frac1d\,T(s+d-1,s).$$ \rm

\bigskip
\noindent
\bf General Triangles. 
\rm 
The above method of proof extends to a much wider class of triangle of numbers generated by the family
$$\{P(x)^nQ(x): n\in\Bbb{N}\}$$
where the polynomial $P(x)$ is palindromic. For the sake of simplicity we will take
$Q(x)=1-x^2$.

\smallskip
\noindent
Fix $d\in\Bbb{N}$ even. Consider for instance the sequence $A(n,k)$ defined by the recurrence
$$A(n,k)=a_0A(n-1,k)+a_1A(n-1,k-1)+\cdots+a_dA(n-1,k-d)$$
satisfying some initial conditions and where $a_j=a_{d-j}$ for $j\in\{0,1,\dots,d\}$ (palindromic coefficients). As before, extend the definition of $A(n,k)$ as skew-symmetric. If we take $P(x)=\sum_{j=0}^da_jx^j$ and $Q(x)=1-x^2$ then
$$\sum_{k\geq0}A(n,k)\,x^k=P(x)^nQ(x).$$
Once more, The Snake Oil method delivers the argument almost verbatim:
$$\align \sum_{k=0}^{dn/2}A(n,k)A(n,k+1)&=\frac12\sum_{k=0}^{dn+2}A(n,k)A(n,k+1) \\
&=\frac12\sum_{k=0}^{dn+2}A(n,k)\,\cdot\text{CT}\left(\frac{P(x)^n(1-x^2)}{x^{k+1}}\right) \\
&=\frac12\,\text{CT}\left[\frac{P(x)^n(1-x^2)}{x}\sum_{k=0}^{dn+2}A(n,k)x^{-k}\right] \\
&=\frac12\,\text{CT}\left[\frac{P(x)^n(1-x^2)}{x}P\left(1/x\right)^n\left(1-\frac1{x^2}\right)\right]  \\
&=\frac12\,\text{CT}\left[\frac{P(x)^{2n}(1-x^2)}{x^{dn+1}}\left(1-\frac1{x^2}\right)\right] \\
&=\frac12\,\text{CT}\left[\frac{P(x)^{2n}(1-x^2)}{x^{dn+1}}\right] -
\frac12\,\text{CT}\left[\frac{P(x)^{2n}(1-x^2)}{x^{dn+3}}\right]  \\
&=\frac12\, A(2n,dn+1)-\frac12\, A(2n,dn+3) \\
&=\frac12\, A(2n,dn-1); \endalign$$
where the last equality is due to $A(2n,dn+1)=0$ and $A(2n,dn+3)=-A(2n,dn-1)$.

\enddocument